\documentclass{amsart}[11pt]
\usepackage{amssymb}

\newtheorem{thm}{Theorem}[section]

\newtheorem{lem}[thm]{Lemma}
\newtheorem{prop}[thm]{Proposition}
\theoremstyle{definition}

\theoremstyle{remark}

\numberwithin{equation}{section}

\newcommand{\R}{\mathbb R}
\newcommand{\C}{\mathbb C}
\newcommand{\Hy}{\mathbb H}
\newcommand{\Cn}{\mathbb C^n}
\newcommand{\B}{\mathbb B^n}
\newcommand{\ep}{\varepsilon}

\newcommand{\om}{\Omega}
\newcommand{\T}{\tilde \Omega}
\newcommand{\po}{\partial \Omega}
\newcommand{\1}{\Omega_1}
\newcommand{\2}{\Omega_2}
\newcommand{\D}{\Delta }
\newcommand{\de}{\delta }
\newcommand{\dx}{\delta (x)}
\newcommand{\dy}{\delta (f(x))}
\newcommand{\dk}{\delta _k}
\newcommand{\gk}{\gamma _k}
\newcommand{\ga}{\gamma}
\newcommand{\si}{\sigma}

\newcommand{\rt}{\rightarrow}

\newcommand{\s}{\rm strongly \ pseudoconvex }
\newcommand{\sit}{\it strongly \ pseudoconvex }
\newcommand{\tz}{ \tilde z}
\newcommand{\tv}{\tilde v}
\newcommand{\bz}{ z_n}
\newcommand{\tw}{\tilde w}
\newcommand{\bw}{ w_n}
\newcommand{\z}{\zeta}
\newcommand{\al}{\alpha}
\newcommand{\3}{\tilde \om_1^k}
\newcommand{\4}{\tilde \om_2^k}
\newcommand{\Z}{ z_1,z_2}
\newcommand{\q}{\rm Carath\'eodory}
\begin{document}
\title[Holomorphicity of Isometries]{On the holomorphicity of isometries of \\ intrinsic metrics in complex analysis}
\author{Harish Seshadri}
\address{department of mathematics,
Indian Institute of Science, Bangalore 560012, India}
\email{harish@math.iisc.ernet.in}
\author{Kaushal Verma}
\address{department of mathematics,
Indian Institute of Science, Bangalore 560012, India}
\email{kverma@math.iisc.ernet.in}

\begin{abstract}
Let $\1$ and $\2$ be $\s$ domains in $\Cn$ and $f: \1 \rt \2$ an
isometry for the Kobayashi or Carath\'eodory metrics. Suppose that
$f$ extends as a $C^1$ map to $ \bar \om_1$. We then prove that
$f|_{\partial \1}: \partial \1 \rt
\partial \2$ is a CR or anti-CR diffeomorphism. It follows that $\1$
and $\2$ must be biholomorphic or anti-biholomorphic. \vspace{1mm}

The main tool is a metric version of the Pinchuk rescaling
technique.
\end{abstract}
\maketitle
\tableofcontents
\section{Introduction}
Complex Finsler metrics such as the Carath\'eodory and Kobayashi
~\cite{kob}
 metrics and K\"ahler metrics such as the Bergman and Cheng-Yau
K\"ahler-Einstein metrics ~\cite{cy} have proved to be very useful
in the study functions of several complex variables. Since
biholomorphic mappings are isometries for these metrics, they are
referred to as ``intrinsic".

This work is motivated by the question of whether
(anti)-biholomorphic mappings are the {\it only} isometries for
these metrics, i.e. is any isometry $f: \1 \rt \2$ between two
domains $\1$ and $\2$ in $\Cn$ (on which the appropriate intrinsic
metrics are non-degenerate) holomorphic or anti-holomorphic ?

To be more precise by what we mean by an isometry, let $F_\om$ and
$d_\om$ denote an intrinsic Finsler metric and the induced
distance on a domain $\om$. In this paper by a $C^0$-isometry we
mean a distance-preserving bijection between the metric spaces $(
\1 ,d_{\1} )$ and $( \2 ,d_{\2} )$. For $k \ge 1$, a
$C^k$-isometry is a $C^k$-diffeomorphism $f$ from $\1$ to $\2$
with $f^\ast (F_{\2})=F_{\1}$. A $C^k$-isometry, $k \ge 1$, is a
$C^0$-isometry and if the Finsler metric comes from a smooth
Riemannian metric (as is the case with the Bergman and the
Cheng-Yau metrics), the converse is also true by a classical
theorem of Myers and Steenrod.

The question above makes sense for a large class of domains (for
example bounded domains). However, we confine ourselves to bounded
strongly pseudoconvex domains in this paper.

We note that the question has been answered in the affirmative
(for strongly pseudoconvex domains) for the Bergman and the
K\"ahler-Einstein metrics in ~\cite{gk1}. The proof is
essentially based on the fact that the metric under consideration
is a K\"ahler metric whose holomorphic sectional curvatures tend
to $-1$ as one approaches the boundary of the domain. Note that
the Bergman metric and the K\"ahler-Einstein metric both have
this property.

The case of the Carath\'eodory and the Kobayashi metrics is more
delicate. A technical reason is that these metrics are Finsler,
not Riemannian, and moreover they are just continuous and not
smooth for general $\s$ domains. Despite these issues, the
results in this paper indicate that the answer to the main
question might be in the affirmative.

Before stating our results we remark that all domains we consider
have at least $C^2$-boundaries. Our main theorem asserts that, an
isometry is indeed a holomorphic mapping at ``infinity".
\begin{thm}\label{mai}
Let $f: \1 \rt \2$ be a $C^1$-isometry of two bounded strongly
pseudoconvex domains in $\Cn$ equipped with the Kobayashi metrics.
Suppose that $f$ extends as a $C^1$ map to $ \bar \Omega_1$. Then
$f|_{\partial \1}:
\partial \1 \rt \partial \2$ is a CR or anti-CR diffeomorphism.
Hence $\1$ and $\2$ must be holomorphic or anti-biholomorphic.

A similar statement holds for the Carath\'eodory metric if we
assume that $\po_1$ and $\po_2$ are $C^3$.
\end{thm}

A few comments about $C^1$-extension assumption: any
$C^1$-isometry between strongly pseudoconvex domains $\Omega_1$
and $\Omega_2$ equipped with the Kobayashi or Carath\'eodory
metrics extends to a $C^ {\frac {1}{2}}$ (H\"older continuous with
exponent $\frac {1}{2}$) map of $\overline \Omega_1$ by the
results of  ~\cite{bb}. A key ingredient in the proof of this
result is that strongly pseudoconvex domains equipped with the
Kobayashi or Carath\'eodory metrics are Gromov hyperbolic.

Our $C^1$-extension assumption is a much stronger one. The
original proof of this extension property for biholomorphisms by
C. Fefferman is based on analysis of the Bergman kernel
~\cite{fef}. It would be interesting to prove the extension
property for Kobayashi/Carath\'eodory isometries and hence render
the assumption in Theorem \ref{mai} unnecessary.


We now summarize the ideas behind the proof of Theorem \ref{mai}
The main idea is to use the rescaling technique of Pinchuk
~\cite{pin3} to study the derivative of the isometry at a boundary
point. We construct a sequence of rescalings of the isometry near
a boundary point $p$ and show that this sequence converges to an
(anti)-holomorphic automorphism of the unbounded realization of
the ball in $\Cn$. On the other hand, we observe that the
horizontal components of these rescalings converge to the
horizontal component of the derivative of the isometry at $p$.
Here by a ``horizontal" vector we mean a vector in the maximal
complex subspace of a tangent space of the domain. In fact, we
show that the restriction of the derivative to the horizontal
subspace at $p$ can be related to the values of the holomorphic
automorphism acting on certain points in the ball. These two
facts together are shown to imply the complex-linearity of the
derivative on the horizontal subspace of the tangent space of
$p$. Much of the technical work in the proof is in showing the
convergence of metrics under Pinchuk rescalings. The first
important technical lemma that we need is about the behaviour of
the distance to the boundary under isometries (Lemma \ref{del}).
Here we use the the two-sided estimates for the Kobayashi
distance obtained in ~\cite{bb}.

\section{Preliminaries}
\subsection{The Kobayashi and Carath\'eodory metrics}
Let $\D$ denote the open unit disc in $\C$ and let $\rho(a,b)$
denote the distance between two points $a,b \in \D$ with respect
to the Poinc\'are metric (of constant curvature $-4$).

Let $\om $ be a domain in $\Cn$. The Kobayashi, \q \ and inner-\q
 \ distances on $\om$, denoted by $d^K_{\om}, \ d^C_{\om}$ and
$d^{\tilde C}_{\om}$ respectively, are defined as follows:

Let $z \in \om$ and $v \in T_z \om$ a tangent vector at $z$.
Define the associated infinitesimal Kobayashi and Carath\'eodory
metrics as
$$F_ \om^K(z,v) =  \; \inf \ \{ \frac {1}{\alpha}: \ \alpha >0, \phi \in {\mathcal O}(
\D, \om) \ {\rm with} \ \phi(0)=z, \ \phi'(0)=\alpha v \}$$ and
$$F_ \om^{\tilde C}(z,v) = \; \sup \ \{ df(z)v:  f \in {\mathcal O}(
 \om, \D) \}.$$
respectively. The inner-Carath\'eodory \ length and the Kobayashi
length of a piece-wise $C^1$ curve $\si :[0,1] \rt \om$ are given
by
$$ l^{\tilde C}(\si)=\int_0^1F_\om^{\tilde C}(\si, \si')dt  \ \
{\rm and} \ \ l^{K}(\si)=\int_0^1F_\om^{K}(\si, \si')dt $$
respectively. Finally the Kobayashi and  inner-\q \ distances
between $p, q$ are defined by $$d_\om^{K}(p,q) = \; \inf \ l^{K
}(p,q) \ \ {\rm and} \ \  d_\om^{\tilde C}(p,q) = \; \inf \
l^{\tilde C}(p,q),$$ where the infimums are taken over all
piece-wise $C^1$ curves in $\Omega$ joining $p$ and $q$.

The \q \ distance is defined to be
$$d_\om^{C}(p,q) = \; \sup \ \{ \rho (f(p),f(q): \ f \in {\mathcal
O}(\om , \D) \}.$$

We note the following well-known and easy facts:
\begin{itemize}

\item[{$\bullet$}] If $\om$ is a bounded domain, then $d^K_{\om},
\ d^C_{\om}$ and $d^{\tilde C}_{\om}$ are non-degenerate and the
topology induced by these distances is the Euclidean topology.

\item[{ $\bullet $}] These distance functions are invariant under
biholomorphisms. More generally, holomorphic mappings are distance
non-increasing. The same holds for $F^K_{\om}(z,v)$ and
$F^{\tilde C}_{\om}(z,v)$.

\item[{ $\bullet $}] We always have
$$ d^C_{\om}(p,q) \le \ d_{\om}^{\tilde C}(p,q) \le d^{K}_{\om}(p,q).$$

\item[{ $\bullet$ }] If $\om=\B$, all the distance functions
above coincide and are equal to the distance function of the
Bergman metric $g_0$ on $\B$. Here The Bergman metric is a
complete K\"ahler metric normailzed to have constant holomorpic
sectional curvature $-4$. Also, for $\B$, the infinitesimal
Kobayashi and Carath\'eodory metrics are both equal to the
quadratic form associated to $g_0$.
\end{itemize}
\subsection{Convexity and Pseudoconvexity}
Suppose $\om$ is a bounded domain in $\Cn$, $n \ge 2$, with
$C^2$-smooth boundary. Let $ \rho: \Cn \rt \R$ be a smooth
defining function for $\om$, i.e, $ \rho =0$ on $\po$, $ d \rho
\neq 0$ at any point of $\po$ and $ \rho ^{-1} ( - \infty, 0) =
\om$.

A domain with $C^2$ smooth boundary $\om$ is said to be {\it
strongly convex} if there is a defining function $\rho$ for $\po$
such that the real Hessian of $\rho$ is positive definite as a
bilinear form on $T_p \po$, for every $p \in \po$.

$\om$ is {\it strictly convex} if the interior of the straight
line segment joining any two points in $ \overline \om$ is
contained in $\om$. Note that we do not demand the boundary of
$\om$ be smooth. Strong convexity implies strict convexity.

Let $\om$ be a bounded domain. A holomorphic map $\phi: \D \rt
\om$ is said to be an {\it extremal disc} or a {\it complex
geodesic} for the Kobayashi metric (or distance) if it is
distance preserving, i.e. $d_\om^K( \phi(p), \phi(q))= \rho(p,q)$
for all $p,q \in \D$.

The following fundamental theorem about complex geodesics in
strictly convex domains will be repeatedly used in Section 3 of
this paper:

\begin{thm} \label{le1}({\bf L. Lempert} ~\cite{lem})
Let $\om$ be a bounded strictly convex domain in $\Cn$.
\vspace{2mm}

(1) Given $p \in \om$ and $\ v \in \Cn$, there exists a complex
geodesic $\phi$ with $\phi (0)=p$ and $\phi '(0)=v$ (or $d \phi
 \ (T_0 \D) = P_v$, where $P_v \subset T_p \om$ is the real -two
plane generated by the complex vector $v$) .

$\phi$ also preserves the infinitesimal metric, i.e., $F^K_\om(
\phi (q) ; d \phi (w)) = F_\D(q;w)$ for all $w \in T_q \D$.
\vspace{2mm}

(2) Given $p$ and $q$ in $\om$, there exists a complex geodesic
$\phi$ whose image contains $p$ and $q$. \vspace{2mm}

(3) The Kobayashi, Carath\'eodory and inner-Carath\'eodory
 distances coincide on $\om$. Also, the Kobayashi and Carth\'eodory
 infinitesimal metrics coincide on $\om$.
\end{thm}



The {\it Levi form} of the defining function $\rho$ at $p \in
\Cn$ is defined by
$$L_p(v)= \sum_{i,j=1}^n \frac { \partial^2 \rho}{\partial z_i
\partial \bar z_j}(p)v_i \overline v_j \ \ {\rm for} \ \ v=(v_1,..,v_n)
\in \Cn.$$

For $p \in \po$, the maximal complex subspace of the tangent
space $T_p \po$ is denoted by $H_p (\po)$ and called the {\it
horizontal subspace} at $p$. By definition, $\om$ is {\it
strongly pseudoconvex} if $L_p$ is positive definite on $H_p(\po)$
for all $p \in \po$. It can be checked that strong convexity
implies strong pseudoconvexity.

For a strongly pseudoconvex domain, the {\it Carnot-Carath\'eodory
metric} on $\po$ is defined as follows. A piecewise $C^1$ curve
$\alpha:[0,1] \rt \po$ is called {\it horizontal} if $\dot \alpha
(t) \in H_{\alpha(t)}(\po)$ wherever $\dot \alpha (t)$ exists. The
strong pseudoconvexity of $\om$ implies that $\po$ is connected
and, in fact, any two points can be connected by a horizontal
curve. The {\it Levi-length} of a curve $\alpha$ is defined by
$l(\alpha) = \int_0^1L_{\alpha(t)}(\dot \alpha (t))^{\frac
{1}{2}}dt$. Finally the {\it Carnot-Carath\'eodory metric} is
defined, for any $p,q \in \po$, by
$$d_H(p,q) = \; \inf_{\alpha} \; l(\alpha),$$
where the infimum is taken over horizontal curves $\alpha:[0,1]
\rt \po$ with $\alpha(0)=p$ and $\alpha(1)=q$.

\subsection{Notation} \hfill

$\circ$ \ \  $\D := \{ z \in \C : \vert z \vert <1 \}$, $\D_r :=
\{z \in \C : \vert z \vert < r \}$.

$\circ$ \ \ $\rho = $ distance function on $\D$ of the Poinc\'are
metric of curvature $-4$.

$\circ$ \ \ For $n \ge 2$, \ $\B:= \{ z \in \Cn : \vert z \vert
<1 \}$ and $B_a(r)= \{z \in \Cn : \vert z -a \vert < r \}$.

$\circ$ \ \ $\Sigma_z = \{ z=(z_1,..,z_n) \in \Cn: \ 2 \ {\rm Re}
\ z_n + \vert z_1 \vert^2+..+\vert z_{n-1} \vert^2 \ <0 \}=$ the

\ \ \ \ unbounded realization of the ball ${\mathbb B}^n$= \
Siegel domain.

$\circ$ \ \ $H_p (\po) \subset T_p \po$ denotes the horizontal
subspace of $T_p \po$, $ p \in \po$.

$\circ$  \ \ Given $p \in \po$, for any $v \in \Cn$, \ $v=
v_H+v_N$ corresponds to $\Cn =$

\ \ \ \ $H_p (\po) \oplus H_p (\po)^{\perp}$.

$\circ$ \ \ $z=( \tz, z_n)$ corresponds to $\Cn= \C ^{n-1} \times
\C$.

$\circ$ \ \ $\dx =d(x,\po)$ denotes Euclidean distance of $x \in
\om$ to $\po$.







$\circ$ \ \ $d^K_{\om}, \ d^C_{\om}$ and $d^{\tilde C}_{\om}$
denote the Kobayashi, \q \ and inner-Carath\'eodory

\ \ \ \ distances on $\om$.

$\circ$ \ \ $F_\om^K$ and $F_\om^C$ denote the Kobayashi and
Carath\'eodory infinitesimal

\ \ \ \ metrics on $\om$.

$\circ$ \ \ If $f: \1 \rt \2$ is a smooth map between domains
$\1$ and $\2$ in $\Cn$,

\ \ \ \ then $df_p: \R^{2n} \rt \R^{2n}$ denotes the derivative
at $p \in \1$. \vspace{3mm}

\noindent Finally, the letters $C$ or $c$ will be used to denote
an arbitrary constant throughout this article and which is subject
to change, even within the limits of a given line, unless
otherwise stated.

\subsection{An estimate for the distance to the boundary}
We prove that $C^0$-isometries approximately preserve the
distance to the boundary. This is needed for the convergence of
Pinchuk rescalings in Section 3. For a domain $\om$ and a point
$x \in \om$, $\dx$ denotes the Euclidean distance $\dx =d(x,
\po)$. Our proof uses the results and notations of ~\cite{bb} in
a crucial way and we refer the reader to it for further details.

We note that in the lemma below we do not need to assume that the
isometry has a $C^1$ extension to the closure of the domain.
\begin{lem}\label{del}
Let $\1$ and $\2$ be strongly pseudoconvex domains in $\Cn$ and
$f: \1 \rt \2$ a $C^0$ isometry of the Kobayashi on $\1$ and
$\2$. There exist positive constants $A$ and $B$ such that
$$B \ \dx \le \ \dy \ \le A \ \dx$$ for all $x \in \1$.
A similar statement holds for an isometry of the Carath\'eodory
distance if we assume that $\po_1$ and $\po_2$ are $C^3$.
\end{lem}
\begin{proof}
Since $\Omega$ has a $C^2$ boundary, given $x \in \Omega$
sufficiently close to the boundary, there exists a unique point
$\pi(x) \in \partial \Omega$ such that $\vert x - \pi(x) \vert=
\dx $. Extend the domain of $\pi$ to be all of $ \Omega$. Such an
extension is not uniquely defined but any extension will work for
our purposes.

Following ~\cite{bb}, define for any $\s$ $\Omega$, the function
$g: \Omega \times \Omega \rt \R$ by
$$g (x,y) = 2 \log \biggl [ \frac {d_H( \pi (x), \pi(y)) \ + \
{\rm max} \{ h(x),h(y) \} }{\sqrt {h(x)h(y)}} \biggr ],$$ where
$h(x) = \dx^{\frac {1}{2}}$ and $d_H$ is the
Carnot-Carath\'eodory metric on $\partial \Omega$.

The Box-Ball estimate (Proposition 3.1 of ~\cite{bb}) implies
that the topology induced by $d_H$ on $\partial \Omega$ agrees
with the Euclidean topology. Hence, $(\partial \Omega, d_H)$ is
compact and, in particular, has finite diameter, say $D$. This
implies that
\begin{equation}
2 \log \sqrt {\frac {h(y)}{h(x)}} \le g(x,y) \le 2 \log \Big(
\frac {D+S}{ \sqrt {h(x)h(y)}} \Big) \notag
\end{equation}
where we have used max$\{h(x),h(y)\} \ge h(y)$ in the first
inequality and where $S = \sup_{x \in \Omega} h(x)$. Hence
\begin{equation}\label{g}
{\frac {h(y)}{h(x)}} \le e^{g(x,y)} \le \frac {E}{h(x)h(y)}
\end{equation}

Now we consider the functions $g_1$ and $g_2$ associated to $\1$
and $\2$.  By Corollary 1.3 of ~\cite{bb}, there exists a
constant $C_1$ such that
\begin{equation}\label{g1}
g_1(x,y) - C_1\le d^K_{\1}(x,y) \le g_1(x,y) + C_1
\end{equation}
for all $x$, $y$ in $\1$.

According to ~\cite{bb1}, such an estimate holds for the
inner-Carath\'eodory distance as well, if one assumes
$C^3$-regularity of the boundaries.

Combining (\ref{g1}) and (\ref{g}) gives
$$ A_1{\frac {h(y)}{h(x)}} \le e^{d^K_{\1}(x,y)} \le \frac {B_1}{h(x)h(y)}.$$
A similar inequality holds on $\2$ (with $A_2$, etc). Fixing $y
\in \1$, using $d^K_{\1}(x,y)=d^K_{\2}(f(x),f(y))$, and comparing
the inequalities on $\1$ and $\2$, we get the required estimates.
The proof for the inner-Carath\'eodory distance is the same.
\end{proof}
An immediate corollary of Lemma \ref{del} is that for
$C^1$-isometries which have $C^1$-extensions, the derivative of
the boundary map preserves the horizontal distribution of $T$.
Note that necessarily $f (\po_1) \subset \po_2$, by Lemma
\ref{del}.
\begin{lem}\label{hor}
Let $f: \1 \rt \2$ be a $C^1$-isometry of strongly pseudoconvex
domains equipped with the Kobayashi metric. If $f$ extends to a
$C^1$-map of $\overline \Omega_1$, then
$$df_p \ ( H_p(\partial \Omega_1) ) \subset  H_{f(p)} \
(\partial \Omega_2),$$ for any $p \in \po$. This holds for an
isometry of the Carath\'eodory metric as well if we assume that
$\po_1$ and $\po_2$ are $C^3$.
\end{lem}
\begin{proof}
By ~\cite{ma}, there exists $\de_0 >0$ such that for any $x \in
\1$ with $ \dx \le \de_0$ and for all $v=v_H+v_N \in \Cn$ (where
this decomposition is taken at $\pi (x)$), we have
\begin{equation}\label{ma}
\frac {\vert v_N \vert ^2} {8 \dx ^2} \ + \ \frac {L_{ \pi (x)}(
v_H)}{4 \dx}  \le \Big( F_{\1}^K(x, v) \Big)^2 \le  \frac {\vert
v_N \vert ^2} {2 \dx ^2} \ + \ 4 \frac {L_{\pi (x)}( v_H)}{ \dx}
\end{equation}
One has similar estimates for $df_x(v)= df_x(v)_H + df_x(v)_N$ .
Since $f$ is an isometry we have $F_{\Omega_1}^K(x, v) =
F_{\Omega_2}^K(f(x), df_x(v))$. Now assume that $v \in H_p (\po_1)
$, i.e., $v = v_H$. Comparing the estimates (corresponding to
(\ref{ma})) for $v$ and $df_x(v)$, we get
\begin{equation}\label{bbb}
 \frac {\vert df_x(v)_N \vert
^2} {8 \dy ^2} \ + \ \frac {L_{ \pi (f(x))}( df_x(v)_H)}{4 \dy}
\le 4 \frac {L_{\pi (x)}( v_H)}{ \dx}
\end{equation}
We can assume that $L_{\pi(x)}( w) \le c \vert w \vert ^2$ for all
$w \in H_q (\po_1), \ q \in \po_1$. Combining this with Lemma
\ref{del} and (\ref{bbb}), we get
$$ \vert df_x(v)_N \vert   \le  C \dx \vert v \vert$$
for some constant $C$. Letting $x \rt p$ and using the continuity
of $df$ we obtain $df_p(v)_N=0$.
\end{proof}
\vspace{2mm}

\section{A metric version of Pinchuk rescaling }
Throughout this section, we will assume that the boundary of the
domain under consideration is $C^3$ when dealing with the
Carth\'eodory distance. Otherwise we assume that $\po$ is $C^2$.

Let $p \in \po$, and fix a sequence $ \{ p_n \}$ in $\om$
converging to $p$. It has been shown in Lemma \ref{del} that
\begin{equation}\label{sim}
d(p_n, \po_1) \sim d (f(p_n), \po_2).
\end{equation}
In particular $\{ f(p_n) \}$ will cluster only on $\po_2$. By
passing to a subsequence we can assume that $q_n=f(p_n) \rt q \in
\po_2$. Fix a defining function $\rho$ for $\po_1$ that is
strongly plurisubharmonic and of class $C^2$ in some neighbourhood
of $\overline  \Omega_1$. Similarly let $\rho'$ be such a function
for $\2$. The following lemma in ~\cite{pin2} will be vital for
what follows.
\begin{lem}\label{vit}
Let $\om$ be a $\sit$ domain, $\rho$ a defining function for
$\po$, and $p \in \po$. Then there exists a neighbourhood $U$ of
$p$ and a family of biholomorphic mappings $h_\z: \Cn \rt \Cn$
depending continuously on $\z \in \po \cap U$ that satisfy the
following:

 (i) \ $h_\z(\z)=0$

 (ii) The defining
function $\rho_\z= \rho \circ h_\z^{-1}$ of the domain
$\om^\z:=h_\z(\om)$ has the form
$$ \rho_\z(z) \ = \ 2  Re(z_n+K_\z(z))+H_\z(z)+\al_\z(z)$$
where $K_\z(z)= \sum^n_{i,j=1}a_{ij}(\z)z_iz_j$, $H_\z(z)=
\sum^n_{i,j=1}a_{i \bar j}(\z)z_i \bar z_j$  and $
\al_\z(z)=o(\vert z \vert ^2)$ with $K_\z( \tz,0) \equiv 0$ and
$H_\z( \tz,0) \equiv \vert \tz \vert^2$

 (iii) The mapping $h_\z$
takes the real normal to $\po$ at $\z$ to the real normal $\{ \tz
=y_n=0 \}$ to $\po^\z$ at the origin.
\end{lem}
To apply this lemma select $\z_k \in \po_1$, closest to $p_k$ and
$w_k \in \po_2$ closest to $q_k=f(p_k)$. For $k$ large, the
choice of $\z_k$ and $w_k$ is unique since $\po_1$ and $\po_2$
are sufficiently smooth. Moreover $\z_k \rt p$ and $w_k \rt q$.
Let $h_k:= h_{\z_k}$ and $g_k:=g_{w_k}$ be the biholomorphic
mappings provided by the lemma above. Let
$$ \1^k:=h_k(\1), \ \ \ \2^k:= g_k( \2) \ \ {\rm and} \ \ f_k:=g_k \circ f \circ h_k^{-1}: \1^k \rt \2^k.$$
Note that $f_k$ is also an isometry for the Kobayashi distance on
$\1^k$ and $\2^k$.

Let $T_k: \Cn \rt \Cn$ be the anisotropic dilatation map given by
$$T_k(\tz,\bz)= (\frac {1}{\sqrt \dk}\tz, \frac {1}{\dk}\bz)$$
and let $$\tilde \om_1^k:=T_k (\om_1^k), \ \ \ \tilde
\om_2^k:=T_k (\om_2^k) \ \ {\rm and} \ \ \Phi_k := T_k \circ f_k
\circ T_k^{-1}: \  \3 \rt \4 .$$ Again $\Phi_k$ is an isometry.
Let us note that the explicit expression for $\Phi_k$ is
$$\Phi_k(z)=\Biggl ( \frac {1}{\sqrt \dk}\tilde f_k(\sqrt \dk  \tz, \dk
\bz), \frac {1}{\dk} \hat {f_k}(\sqrt \dk \tz, \dk \bz) \Biggl)
.$$

For notational convenience, let us denote the compositions of the
rotations and the scalings by
\begin{equation}\label{ba}
H_k:= T_k \circ h_k \ \ {\rm and} \ \ I_k:= T_k \circ g_k \ \
\Rightarrow  \ \ \Phi_k = I_k \circ f \circ H_k.
\end{equation}

 Note that the defining functions for $\tilde \om_1^k$ and
$\tilde \om_2^k$ are given by
 $$\tilde \rho_k(z)=  \frac {1}{\sqrt \dk} \rho _k(\sqrt \dk \tz, \dk \bz),
 \ \ \
\tilde \rho' _k(w)=  \frac {1}{\sqrt \dk} \rho _k(\sqrt \dk \tw,
\dk \bw)$$ respectively.

The family of functions $\{ h_k \}$ converges uniformly on compact
subsets of $\Cn$ to the identity mapping, as do their inverses $
h_k^{-1} $. Thus it follows that for $k>>1$
\begin{equation}\label{sim1}
\frac {1}{B}  \ \le \ \frac {d (  h_k(z), \   \po_1^k )}{ d(z,
 \ \po_1)} \ \le \ B \ \ {\rm and } \ \
 \frac {1}{B}  \ \le \
\frac {d(g_k(w), \ \po_2^k)}{ d(w, \ \po_2)} \ \le \ B
\end{equation}
for some constant $B$ independent of $z$ and $k$.

Combining (\ref{sim}) and (\ref{sim1}) shows that for $k>>1$
\begin{equation}\label{sim3}
 \frac {1}{c} \ d(z, \ \po_1^k)
\ \le \ d(f_k(z),  \ \po_2^k) \ \le \  c \ d(z, \ \po_1^k)
\end{equation}
 where $c$ is independent of $k$ ( for $k>>1$) and $z \in
\1^k$. Moreover, since $\rho$ and $\rho'$ (and hence $\rho_k:=
\rho_{\z_k}=\rho \circ h_k^{-1}$ and $\rho'_k:=\rho'_{w_k}=\rho'
\circ g_k^{-1}$) are smooth, it follows that there exists a
uniform constant $c >0$ such that
\begin{equation}\label{sim4}
\frac {1}{c}  \ \le \ \frac {d(z, \ \po_1^k)}{ \vert \rho_k(z)
\vert} \ \le \ c \ \ {\rm and} \ \
\frac {1}{c}  \ \le \ \frac {d(w, \ \po_2^k)}{\vert \rho'_k(w)
\vert}  \ \le \ c
\end{equation}
for $k>>1$ and $z \in \1^k, \ w \in \2^k$. Let $\dk= d(h_k(p_k),
\ \po_1^k)$ and $ \gk = d(g_k(q_k), \ \po_2^k)$.

Two observations can be made at this stage: first, for $k>>1$,
$h_k(p_k)=(0, - \dk), \ g_k(q_k)=(0,- \gk)$ and $f_k(0, -
\dk)=(0, - \gk)$ as follows from Lemma \ref{vit}, and secondly
(\ref{sim3}) shows that
\begin{equation}\label{sim6}
\frac {1}{c} \ \dk \ \le \ \gk \ \le \ c \ \dk
\end{equation}
for some $c>0$.

It has been shown in ~\cite{pin2} that the sequence of domains $\{
\tilde \om_1^k \}$ converges to the unbounded realization of the
unit ball, namely to
$$\Sigma_z= \{ z \in \Cn : \ 2{\rm Re} \ \bz + \vert \tz \vert ^2 <0 \}.$$
The convergence is in the sense of Hausdorff convergence of sets.
Similarly $\{ \tilde \om_2^k \}$ will converge to $ \Sigma_w$,
the unbounded realization of the ball in $w$ coordinates.
\begin{prop}
Let $\1, \2$ be smoothly bounded strongly pseudoconvex domains in
$\Cn$. Suppose that $f: \1 \rt \2$ is a $C^0$-isometry with
respect to the Kobayashi (resp. Carath\'eodory) distances on $\1$
and $\2$. Define the sequence of domains $\tilde \om_1^k, \ \tilde
\om_2^k$ and mappings $\Phi_k: \tilde \om_1^k \rt \tilde \om_2^k$.
Then there exists a subsequence of $\{ \Phi_k \}$ that converges
uniformly on compact subsets of $\Sigma_z$ to a continuous mapping
$\Phi: \Sigma_z \rt \Cn$.
\end{prop}
\begin{proof}
The case when $\Phi$ is an isometry with respect to Kobayashi
distances will be dealt with first. By construction
$\Phi_k(0,-1)=(0,- \frac{\gk}{\dk})$ and (\ref{sim6}) shows that
$\{ \Phi_k(0,-1) \}$ is bounded. The domain $\Sigma_z$ can be
exhausted by an increasing union $ \{ S_i \}$ of relatively
compact convex domains each containing $(0,-1)$. Fix a pair
$S_{i_0} \subset \subset S_{i_0+1}$ and write $S_1 =S_{i_0}$ and
$S_2 =S_{i_0+1}$ for brevity. Since $\tilde \om_1^k$ converges to
$\Sigma_z$ it follows that $S_1 \subset \subset S_2 \subset
\subset \tilde \om_1^k$ for all $k >> 1$. It will suffice to show
that $\{ \Phi_k \}$ restricted to $S_1$ is uniformly bounded and
equicontinuous. Fix $s_1,s_2$ in $S_1$. The following
inequalities hold for large $k$:
\begin{equation}\label{koin}
d_{\4}^K(\Phi_k(s_1),\Phi_k(s_2))= d_{\3}^K(s_1,s_2) \le
d_{S_2}^K(s_1,s_2) \le c \vert s_1 - s_2 \vert
\end{equation}
for $c>0$ independent of $k$. Indeed the equality holds for all
$k$ since $\Phi_k$ is an isometry and the inequalities are a
result of the following observations: first, the inclusion $S_2
\hookrightarrow \3$ is distance decreasing for the Kobayashi
distance and second, since $S_2$ is convex, the infinitesimal
Kobayashi metric $F_{S_2}^K(z,v)$ satisfies
\begin{equation}\label{inko}
F_{S_2}^K(z,v) \le \frac {\vert v \vert}{\de_v(z)}
\end{equation}
where $z \in S_2$, $v$ is a tangent vector at $z$ and $\de_v(z)$
is the distance of $z$ to $\partial S_2$ in the direction along
$v$. Now joining $s_1$ and $s_2$ by a straight line path $\gamma
(t)$ and integrating (\ref{inko}) along $\gamma (t)$ yields the
last inequality in (\ref{koin}).

To estimate $d_{\4}^K(\Phi_k(s_1),\Phi_k(s_2))$ note from
(\ref{ba}) that
\begin{align}\label{bor1}
d_{\4}^K(\Phi_k(s_1),\Phi_k(s_2))
= d_{\2}^K(f \circ H_k^{-1}(s_1), \ f \circ H_k^{-1} (s_2))
\end{align}
since $I_k$ is an isometry. Since $f$ is continuous at $p \in
\po_1$, choose neighborhoods $U_1, \ U_2$ of $p, \ q = f(p)$
respectively so that $f(U_1 \cap \1) \subset U_2 \cap \2$. Note
that $p_n, $ and $\z_n$ as chosen earlier lie in $ U_1 \cap \1$
eventually. For $k$ large, $H_k^{-1}(S_1) \subset U_1 \cap \1$
and hence both $f \circ H_k^{-1}(s_1)$ and $f \circ
H_k^{-1}(s_2)$ lie in $U_2 \cap \2$. It is well known that the
Kobayashi distance can be localized near strongly pseudoconvex
points in the sense that for every choice of $U_2$, there is a
smaller neighborhood $p \in U_3$, $U_3$ relatively compact in
$U_2$, and $c > 0$ such that
\begin{equation}\label{bor2}
c d_{U_2 \cap \2}^K(x,y) \le d_{\2}^K(x,y)
\end{equation}
for all $x,y \in U_3 \cap \2$. We apply this to $x=f \circ
H_k^{-1}(s_1)$ and $y= f \circ H_k^{-1} (s_2)$, both of which
belong to $U_3 \cap \Omega_2$ for large $k$, by shrinking $U_1$
if necessary.

Moreover, thanks to the strong pseudoconvexity of $\po_2$ near
$q$, it is possible to choose $U_2$ small enough so that for $k
>> 1$,
$$g_k(U_2 \cap \2) \subset \{ w \in \Cn: \ \vert  w_n + R \vert^2 +
\vert \tilde w \vert^2 <R^2 \} \subset \tilde \om$$ where $\om :=
\{ w \in \Cn : 2 \;R \;({\rm Re} w_n) < - \vert \tilde w \vert ^2
\}$ for some fixed $R > 1$.

Note that $\tilde \om$ is invariant under the dilatation $T_k$
for all $k$ and moreover $\tilde \om$ is biholomorphic to $\B$.
Thus, $T_k \circ g_k(U_2 \cap \2) \subset \tilde \om$ and hence
$\Phi_k(s_1), \Phi_k(s_2)$ both lie in $ \tilde \om$ for $k$
large. From (\ref{bor1}) and (\ref{bor2}) it follows that
\begin{equation}\label{bor3}
c d_{\tilde \om}^K(\Phi_k(s_1),\Phi_k(s_2)) \le
d_{\4}^K(\Phi_k(s_1),\Phi_k(s_2))
\end{equation}
for $k$ large. Combining (\ref{koin}) and (\ref{bor3}) gives
\begin{equation}
d_{\tilde \om}^K(\Phi_k(s_1),\Phi_k(s_2)) \le c \vert s_1 - s_2
\vert
\end{equation}
for $s_1,s_2 \in S_1$ and $k >> 1$.

Let $\psi: \T \rt \B$ be a biholomorphic mapping. To show that
$\{ \Phi_k(S_1) \}$ is uniformly bounded, choose $s_1 \in S_1$
arbitrarily and $s_2 = (0,-1)$. Then (\ref{bor3}) shows that
$$d_{\tilde \om}^K(\Phi_k(s_1),\Phi_k(0,-1)) \le c \vert s_1 - s_2 \vert <
\infty .$$ Since $ \{ \Phi_k(0,-1) \}$ is bounded and $ \B$ (and
hence $\T$) is complete in the Kobayashi distance, it follows that
$\{ \Phi_k (s_1) \}$ is bounded.

To show that $\{ \Phi_k \}$ restricted to $S_1$ is equicontinuous
observe that the Kobayashi distance in $\B$ between $\psi \circ
\Phi_k(s_1)$ and $\psi \circ \Phi_k(s_2)$ equals $d_{ \tilde
\om}^K(\Phi_k(s_1),\Phi_k(s_2)) \le c \vert s_1 - s_2 \vert$.
Using the explicit formula for the Kobayashi distance between two
points in $\B$, this gives
$$ \Bigl \vert  \frac {\psi \circ \Phi_k(s_1) - \psi \circ \Phi_k(s_2)}{ 1 - \overline {\psi \circ \Phi_k(s_1)} \ \psi \circ
\Phi_k(s_2)} \Bigl \vert \ \le \ \frac {{\rm exp} \ (2 c \vert s_1
-s_2 \vert )-1}{{\rm exp} \ (2 c \vert s_1 -s_2 \vert )+1}.$$
Since $\{ \Phi_k(S_1) \}$ is relatively compact in $\T$ for
$k>>1$, it follows that so is $\{ \psi \circ \Phi_k(S_1) \}$ in
$\B$. Let $\{ \psi \circ \Phi_k(S_1) \} \subset G \subset \subset
\B$. Hence
$$ \vert 1 - \overline { \psi \circ \Phi_k(s_1)} \psi \circ \Phi_k(s_2) \vert  \ge c >0$$
for $k$ large and this shows that
$$ \vert \psi \circ \Phi_k(s_1) - \psi \circ \Phi_k(s_2) \vert \le \frac { {\rm exp} \ (c \vert s_1 - s_2 \vert
)-1}{{\rm exp} \ (c \vert s_1 - s_2 \vert )+1} \le c \vert s_1 -
s_2 \vert $$
 This shows that $ \{ \Phi_k \}$ is equicontinuous on $S_1$ and hence
 there is a subsequence of $\{ \Phi_k \}$ that converges uniformly on
 compact subsets of $\Sigma_z$ to a continuous mapping $\Phi:
 \Sigma_z \rt \Cn$.

It may be observed that the same proof works when $f: \1 \rt \2$
is an isometry in the Carath\'eodory distance on the domains.
Indeed, the process of defining the scaling does not depend on
the distance function used. Moreover the Carath\'eodory distance
enjoys the same functorial properties as the Kobayashi distance
and even the quantitative bounds used in (\ref{inko}) and
(\ref{bor2}) remain the same. Hence the same proof works verbatim
for the Carath\'eodory distance.
\end{proof}

\begin{prop}\label{con}
Let $f: \1 \rt \2$ be an isometry in the Kobayashi distance on
smoothly bounded strongly pseudoconvex domains $\1, \2$ in $\Cn$.
Then the limit map map $\Phi: \Sigma_z \rt \Cn$ constructed above
satisfies:\\
(i) $\Phi (\Sigma_z) \subset \Sigma_w$\\
(ii) $\Phi : \Sigma_z \rt \Sigma_w$ is a $C^0$-isometry for the
Kobayashi distance.

The same conclusions hold when $f: \1 \rt \2$ is an isometry in
the Carath\'eodory distance on the domains. In particular $\Phi:
\Sigma_z \rt \Sigma_w$ is an isometry in the Carath\'eodory
distance (which equals the Kobayashi distance) on $\Sigma_z$ and
$\Sigma_w$.
\end{prop}
\begin{proof}
Let $\Phi_k: \3 \rt \4$ be the sequence of scaled mappings as
before. Without loss of generality we assume that $\Phi_k \rt
\Phi$ uniformly on compact subsets of $\Sigma_z$. The defining
equations for $\4$ and $\3$ are respectively given by
$$\tilde \rho' _k(w)=  \frac {1}{\sqrt \dk} \rho _k(\sqrt \dk \tw, \dk
\bw) \ , \ \ \tilde \rho _k(z)=  \frac {1}{\sqrt \dk} \rho
_k(\sqrt \dk \tz, \dk \bz) .$$

It is shown in ~\cite{pin2} that these equations simplify as
$$ \tilde \rho' _k(w)=2 \;{\rm Re} w_n + \vert \tilde w \vert^2 + \tilde
B_k(w) \ , \ \ \tilde \rho _k(w) = 2 \; {\rm Re} z_n + \vert
\tilde z \vert^2 + \tilde A_k(w)$$ in neighborhoods of the origin
where
$$\vert \tilde B_k(w) \vert \le \vert w \vert^2(c \sqrt \dk + \eta
(\dk \vert w \vert^2)) \ , \ \ \vert \tilde A_k(z) \vert \le
\vert z \vert^2(c \sqrt \dk + \eta (\dk \vert z \vert^2))$$ with
$ \eta(t) = o(1)$ as $t \rt 0$ and $c > 0$ is uniform for all $k$
large.

Fix a compact subset of $\Sigma_z$, say $C$. Then for $k>>1$ and
$z \in C$
\begin{equation}\label{def}
2\; {\rm Re}(\Phi_k)_n(z) + \vert \tilde \Phi_k(z) \vert^2 +
\tilde B_k( \Phi _k(z))<0
\end{equation}
where
$$ \vert \tilde B_k(\Phi_k(z)) \vert \le \vert \Phi_k(z) \vert^2(c \sqrt \dk + \eta
(\dk \vert \Phi_k(z) \vert^2)).$$ By the previous proposition  $\{
\Phi_k \}$ is uniformly bounded on $C$ and hence $ \dk \vert
\Phi_k(z) \vert ^2 \rt 0$ with the result that $\eta (\dk \vert
\Phi_k(z) \vert^2)) \rt 0$ as $k \rt \infty$. Passing to the limit
as $k \rt \infty $ in (\ref{def}) shows that
$$ 2{\rm Re} \Phi_n(z)+ \vert \tilde \Phi (z) \vert^2 \le 0$$
which means exactly that $\Phi(C) \subset \overline \Sigma_w$.
Since $C \subset \Sigma_z$ is arbitrary it follows that $\Phi:
\Sigma_z \rt \overline \Sigma_w$. If $\Phi$ were known to be
holomorphic it would follow at once by the maximum principle that
$\Phi: \Sigma_z \rt \Sigma_w$. However $\Phi$ is known to be just
continuous.

Let $D \subset \Sigma_z$ be the set of all points $z$ such that
$\Phi(z) \in \Sigma_w$. $D$ is non-empty since $(0,-1) \in D$ as
can be seen from (\ref{sim6}) and the fact that $\Phi_k(0,-1) =
(0, - \frac{\gk}{\dk})$. Since $\Phi$ is continuous, $D$ is open
in $\Sigma_z$.\\

{\it Claim:} It suffices to show that
\begin{equation}\label{ag1}
d^K_{\Sigma_z}(z_1,z_2) =d^K_{\Sigma_w}(\Phi(z_1),\Phi(z_2))
\end{equation}
for $z_1,z_2 \in D$. Indeed, if $z_0 \in \partial D \cap
\Sigma_z$, choose a sequence $z_j \in D$ that converges to $z_0$.
If the claim were true, then
\begin{equation}\label{ag2}
d^K_{\Sigma_z}(z_j,(0,-1)) =d^K_{\Sigma_w}(\Phi(z_j),\Phi(0,-1))
\end{equation}
for all $j$. Since $z_0 \in \partial D, \ \Phi(z_j) \rt \partial
\Sigma_w$ and as $\Sigma_w$ is complete in the Kobayashi distance,
the right side in (\ref{ag2}) becomes unbounded. However the left
side remains bounded, again because of the completeness of
$\Sigma_z$. This contradiction would show that $D= \Sigma_z$,
knowing which the claim would prove assertion $(ii)$ as well.

It is already known that
$$d^K_{\3}(z_1,z_2) =d^K_{\4}(\Phi(z_1),\Phi(z_2))$$
for $k>>1$. To prove the claim it suffices to take limits on both
sides in the equality above. This is an issue of the stability of
the Kobayashi distance, to understand which we need to study the
behaviour of the infinitesimal Kobayashi metric $\Phi_{\3}(z,v)$
as $k \rt \infty$. To do this, we will use ideas from ~\cite{gk2}.
Once this is done, an integration argument will yield information
about the global metric $K_{\3}$.\\

{\it Step 1:} It will be shown that
\begin{equation}\label{dif}
\lim_{k \rt \infty} F^K_{\3}(a,v)=F^K_{\Sigma_z}(a,v)
\end{equation}
for $(a,v) \in \Sigma_z \times \Cn$. Moreover, the convergence is
uniform on compact subsets of $\Sigma_z \times \Cn$.

Let $S \subset \Sigma_z$ and $C \subset \Cn$ be compact and
suppose that the desired convergence does not occur. Then there
are points $a_k \in S$ converging to $a \in S$ and vectors $v_k
\in G$ converging to $v \in G$ such that
$$0 < \ep _0 <  \vert F^K_{\3} (a_j, v_j) -F^K_{\Sigma_z}(a_j,v_j)
\vert $$ for $j$ large. This inequality holds for a subsequence
only, which will again be denoted by the same symbols. Further,
since the infinitesimal metric is homogeneous of degree one in the
vector variable , we can assume that $\vert v_j \vert=1$ for all
$j$. It was proved in ~\cite{gk2} that $F^K_{\Sigma_z}$ is jointly
continuous in $(z,v)$. This was a consequence of the fact that
$\Sigma_z $ is taut. Thus
$$0 <  {\ep _0} / 2 <  \vert F^K_{\3} (a_j, v_j) -F^K_{\Sigma_z}(a,v)
\vert .$$ The tautness of $\Sigma_{\Phi_z}$ implies, via a normal
families argument, that $0 < F^K_ \Sigma(a,v) < \infty$ and there
exists a holomorphic extremal disc $\phi: \D \rt \Sigma_z$ that
by definition satisfies $\phi (0) =a, \ \phi ' (0)=\mu v$ where
$\mu >0$ and $F^K_ \Sigma (a,v)= \frac {1}{\mu}$.

Fix $\de \in (0,1)$ and define the holomorphic maps $\phi_k: \D
\rt \Cn$ by
$$ \phi_k (\xi) = \phi ((1- \de ) \xi) + (a_k-a) + \mu(1 - \de )
\xi (v_k -v).$$ Observe that the image  $\phi ((1 - \de ) \xi )$
is compact in $\Sigma_z$ and since $a_k \rt a, \ v_k \rt v$ it
follows that $\phi _k : \D \rt \3$ for $k$ large. Also,
$\phi_k(0)= \phi(0) + a_k - a =a_k$ and that $\phi'_k (0)= (1 -
\de ) \phi'(0) +  \mu (1 - \de )(v_k - v) = \mu (1 - \de )(v + v_k
- v)= \mu (1 - \de )v_k$.

By the definition of the infinitesimal metric it follows that
$$F_{\3}(a_k, v_k) \le \frac {1}{\mu (1- \de)}= \frac
{F^K_{\Sigma_z}(a,v)}{1 - \de}$$ for $j >>1$. Letting $\de \rt
0^+$ it follows that
\begin{equation}
\limsup_{ \ k \rt \infty} F ^K _{\3}(a_k,v_k) \ge F ^K
_{\Sigma_z}(a,v).
\end{equation}
Conversely, fix $\ep >0$ arbitrarily small. By definition, there
are holomorphic mappings $\phi_k: \D \rt \3$ satisfying $\phi_k
(0) =a_k, \ \phi_k ' (0)=\mu_k v_k$ where $\mu_k
>0$ and
\begin{equation}\label{sho}
F^K_{\3} (a_k,v_k) \ge \frac {1}{\mu_k} -\ep.
\end{equation}
The sequence $\{ \phi_k \}$ has a subsequence that converges to a
holomorphic mapping $ \phi : \D \rt \Sigma_z$ uniformly on
compact subsets of $\D$. Indeed consider the disc $\D_r$ of
radius $r \in (0,1)$. The mappings $H_k^{-1} \circ \phi_k : \D \rt
\1$ and $H_k^{-1}  \circ \phi_k (0) \rt p \in \po_1$. Fix a ball
$B_p(\de )$ of radius $\de$ around $p$, with $\de$ small enough.
Since $p \in \po_1$ is a plurisubharmonic peak point, Proposition
5.1 in ~\cite{ver} (see ~\cite{ber} also, where this phenomenon
was aptly termed {\it the attraction property of analytic discs})
shows that for the value of $r \in (0,1)$ fixed earlier, there
exists $\eta >0$, independent of $\phi_k$ such that $H_k^{-1}
\circ \phi_k ( \D_r) \subset B_p(\de).$ If $\de$ is small enough,
then there exists $R > 1$ large enough so that
\[
h_k(B_p(\de) \cap \1) \subset \{ z \in \Cn : \vert z_n + R \vert
^2 + \vert \tz \vert^2 < R^2 \} \subset  \Omega.
\]
where (as in Proposition 3.2)
\[
\Omega = \{ z \in \Cn : 2 \; R \; ({\rm Re} \; z_n) < - \vert
\tilde z \vert ^2 \}.
\]

 Again, we note that $\om$ is
invariant under $T_k$ and that $\om  \cong \B$. Hence
$\phi_k(\D_r) \subset \om$ for $k$ large and this exactly means
that
$$ 2\; R \; ({\rm Re} (\phi_k)_n (z)) + \vert \tilde{\phi}_k(z) \vert^2 <
0$$ for $z \in \D_r$.

It follows that $\{ (\phi_k)_n(z) \}$ and hence each component of
$\tilde{\phi}_k(z)$,  forms a normal family on $\D_r$. Since $ r
\in (0,1)$ was arbitrary, the usual diagonal subsequence yields a
holomorphic mapping $\phi: \D \rt \Cn$ or $ \phi \equiv \infty$
on $\D$. But it is not possible that $\phi \equiv \infty$ on $\D$
since $\phi(0) \rightarrow a$.

It remains to show that $ \D \rt \Sigma_z$. For this note that
$\3$ are defined by
\begin{equation}\label{aga}
 \tilde \rho _k(w) =2{\rm Re} \ z_n + \vert \tilde z \vert^2 + \tilde
A_k(w),
\end{equation}
where
$$ \vert \tilde A_k(z) \vert \le \vert z \vert^2(c
\sqrt \dk + \eta (\dk \vert z \vert^2))$$ Thus for $ \z \in \D_r,
\ r \in (0,1)$,
$$ 2 \;R \;({\rm Re} (\phi_k)_n (z) + \vert  {\tilde \phi}_k(z)  \vert^2 + \tilde
A_k( \phi_k(z)) < 0$$ where
$$ \tilde A_k( \phi_k(z)) \le \vert  \phi_k (z) \vert^2(c
\sqrt \dk + \eta (\dk \vert  \phi_k(z) \vert^2))$$ as $k \rt
\infty$. Passing to the limit in (\ref{aga}) shows that $$2 \; R
\;( {\rm Re}  (\phi_k)_n (z)) + \vert  {\tilde \phi}_k(z)  \vert^2
\le 0$$ for $z \in \D_r$, which exactly means that $\phi (\D_r)
\subset \bar \Sigma_z$. Since $r \in (0,1)$ was arbitrary, it
follows that $\phi (\D) \subset \bar \Sigma_z$ and the maximum
principle shows that $\phi (\D) \subset \Sigma_z$.

Note that $\phi (0) = a$ and $\phi '(0) = \lim_{ \ k \rt \infty}
\phi'_k(0) = \lim _{ \ k \rt \infty} \mu_k  v_k = \mu v$ for some
$\mu
>0$. It follows from (\ref{sho}) that
\begin{equation}\label{str}
\liminf _{ \ k \rt \infty} F^K_{\3}(a_k, v_k) \ge
F^K_{\Sigma_z}(a,v) - \ep.
\end{equation}
Combining (\ref{sho}) and (\ref{str}) shows that
$$ \lim  _{ \ k \rt \infty} F^K_{\3}(a_k, v_k) =
F^K_{\Sigma_z}(a,v)$$ which contradicts the assumption made and
proves (\ref{dif}).\\

{\it Step 2:} The goal will now be to integrate (\ref{dif}) to
recover the behaviour of the global metric, i.e. the distance
function.

Let $\ga:[0,1] \rt \Sigma_z$ be a $C^1$ such that $\ga (0) = z_1$
and $\ga (0) = z_2$ and
$$d^K_{\Sigma_k}(z_1,z_2)= \int_0^1 F^K_{\Sigma_z}(\ga, \ga')dt.$$
Since $\ga \subset \subset \Sigma_z, \ \ga \subset \subset \3$
for $k$ large. By Step 1, it follows that
\begin{align}\notag
\int_0^1 F^K_{\3}(\ga, \ga')dt  \le \int_0^1 F^K_{\Sigma_z}(\ga,
\ga')dt + \ep
                              = d^K_{\Sigma_z}(z_1,z_2) + \ep
                             \notag
\end{align}
By definition of $K_{\3}(z_1,z_2)$ it follows that
$$d^K_{\3}(z_1,z_2) \le \int_0^1 F^K_{\3}(\ga, \ga')dt \le
d^K_{\Sigma_z}(z_1,z_2)+ \ep .$$ Thus
\begin{equation}\label{iri}
\limsup_{ \ k \rt \infty} \ d^K_{\3}(z_1,z_2) \ \le \
d^K_{\Sigma_z}(z_1,z_2)
\end{equation}
Conversely since $z_1,z_2 \in D \subset \Sigma_z$, it follows
that $z_1,z_2 \in \3$ for $k$ large. Fix $\ep >0$ and let
$B_p(\eta_1)$ be a small enough neighbourhood of $p \in \po$.
Choose $\eta_2 < \eta _1$ so that
\begin{equation}\label{iri2}
F^K_{\1}(z,v) \le F^K_{B_p(\eta_1) \cap \1} (z,v) \le (1 +
\ep)F^K_{\1}(z,v).
\end{equation}
for $z \in B_p(\eta_2) \cap \1$ and $v$ a tangent vector at $z$.
This is possible by the localization property of the Kobayashi
metric near strongly pseudoconvex points.

If $k$ is sufficiently large, $H_k^{-1}  (z_1)$ and $H_k^{-1}
(z_1)$ both belong to $B_p(\eta_2) \cap \1$. If $\eta_1$ is small
enough, $B_p(\eta_1) \cap \1$ is strictly convex and it follows
from Lempert's work that there exist $m_k > 1$ and holomorphic
mappings
$$ \phi_k: \D_{m_k} \rt B_p(\eta_1) \cap \1,$$
such that $\phi_k(0)= H_k^{-1}(z_1)$, $ \phi_k(1)=H_k^{-1} (z_2)$
and
\begin{equation}\label{iri3}
d^K_{B_p(\eta_1) \cap \1} (H_k^{-1} (z_1), H_k^{-1} (z_2)) =
\rho_{\D_{m_k}}(0, 1) = \int_0^1 F^K_{B_p(\eta_1) \cap \1}(
\phi_k(t), \phi'(t))dt.
\end{equation}
Integrating (\ref{iri2}) and using the fact that $H_k$ are
biholomorphisms and hence Kobayashi isometries, it follows that
$$d^K_{H_k (B_p(\eta_1) \cap \1)} (z_1,z_2) \le (1+ \ep)
d^K_{\3}(z_1,z_2).$$ Hence (\ref{iri3}) shows that
\begin{align}\notag
 \frac {1}{2} \log \Big( \frac {m_k +1}{m_k -1} \Big) =
\rho_{ \D_{m_k}}(0,1)
 &= d^K_{B_p(\eta_1) \cap \1}(H_k ^{-1} (z_1),
 H_k^{-1}(z_2)) \\ \notag  & = d^K_{H_k(B_p(\eta_1)
\cap \1)}(z_1,z_2) \le (1 + \ep) d^K_{\3}(z_1,z_2). \notag
\end{align}

But $$d^K_{\3} (z_1,z_2) \le d^K_{\Sigma_z}(z_1,z_2) + \ep <
\infty $$ and hence $m_k > 1 + \de$ for some uniform $\de > 1$ for
all $k
>>1$. Thus the holomorphic mappings $\sigma_k:= T_k \circ h_k
\circ \phi_k : \D_{1+ \de} \rt T_k \circ h_k (U \cap \1) \subset
\3$ are well-defined and satisfy $\sigma_k(0)= z_1$ and
$\sigma_k(1) = z_2$. Now exactly the same arguments that were used
to establish the lower semi-continuity of the infinitesimal metric
in Step 1 show that $\{ \sigma_k \}$ is a normal family and
$\sigma_k \rt \sigma: \D_{1+ \de} \rt \Sigma_z$ uniformly on
compact subsets of $\de_{1+ \de}$. Again using (\ref{iri2}) and
(\ref{iri3}) we get
\begin{align}\notag
 \int_0^1 F^K_{\3}(\sigma_k, \sigma_k')dt  \le \int_0^1 F^K_{H_k
(B_p(\eta_1) \cap \1)}(\sigma_k, \sigma_k')dt & = d^K_{H_k
(B_p(\eta_1) \cap \1)}(z_1,z_2) \\ \notag & \le (1 +\ep)
d^K_{\3}(z_1,z_2).\notag
\end{align}

 Since $\sigma_k \rt \sigma$ and
$\sigma_k' \rt \sigma'$ uniformly on $[0,1]$, Step 1 shows that
$$\int_{0}^{1} F^K_{\Sigma_z}( \sigma, \sigma ) dt \le \int_0^1 F^K_{\3}(\sigma_k, \sigma_k')dt + \ep \le
d^K_{\Sigma_z}(z_1,z_2)+ C \ep $$ for all large $k$.

It remains to note that since $\sigma (t), \ 0 \le t \le 1$ joins
$z_1,z_2$ it follows that
\begin{equation}\label{bla}
d^K_{\Sigma_z}(z_1,z_2) \le \int_0^1 F^K_{\3}(\sigma_k,
\sigma_k')dt \le d^K_{\3}(z_1,z_2) +C \ep
\end{equation}
Combining (\ref{iri}) and (\ref{bla}), we see that
$$d^K_{\3}(z_1,z_2) \rt d^K_{\Sigma_z}(z_1,z_2)$$
for all $z_1,z_2 \in D$. Exactly the same arguments show that it
is possible to pass to the limit on the right side of (3.14). The
claim made in (\ref{ag1}) follows.

To complete the proof of the proposition for the Kobayashi
distance, it remains to show that $\Phi: \Sigma_z \rt \Sigma_w$ is
surjective. This follows by repeating the argument of the previous
proposition for $f^{-1}: \2 \rt \1$ and considering the scaled
inverses, i.e. $\Psi_k = \Phi_k^{-1}: \4 \rt \3$. This family
will converge to a continuous map $\Psi : \Sigma_w \rt \Cn$
uniformly on compact subsets of $\Sigma_w$. The arguments of this
proposition will then show that $\Psi$ maps $\Sigma_w$ to
$\Sigma_z$. Finally observe that for $w$ in a fixed compact set
$C \subset \Sigma_w$,
\begin{align}\notag
\vert w  - \Phi \circ \Psi (w) \vert & = \vert \Phi_k \circ \Psi_k
(w) - \Phi \circ \Psi (w) \vert  \\ \notag & \le \vert \Phi_k
\circ \Psi_k (w) - \Phi \circ \Psi_k (w) \vert + \vert \Phi \circ
\Psi_k (w) - \Phi \circ \Psi (w) \vert \rt 0 \notag
\end{align}
as $k \rt \infty$. Thus $\Phi \circ \Psi = id = \Psi \circ \Phi$.

We now deal with the case when $f: \1 \rt \2$ is an isometry for
the Carth\'eodory distance on $\1$ and $\2$.

One possibility is to first show that
$$\lim _{k \rt \infty}d^C_{\3}(z_1,z_2)=d^C_{\Sigma_z}(z_1,z_2)$$
for $z_1,z_2 \in \Sigma_z$. Knowing this, the following
inequalities hold:
\begin{align}\notag
d^C_{\Sigma_z}(z_1,z_2) = \lim_{k \rt \infty} \ d^{
C}_{\3}(z_1,z_2) & \le \lim_{k \rt \infty} \ d^{\tilde
C}_{\3}(z_1,z_2) \\ \notag & \le \lim_{k \rt \infty} \
d^K_{\3}(z_1,z_2)= d^K_{\Sigma_z}(z_1,z_2) \notag
\end{align}
Since $d^C_{\Sigma_z}=d^{\tilde C}_{\Sigma_z}=d^K_{\Sigma_z}$, it
would follow that
$$d^{\tilde c}_{\Sigma_z}(z_1,z_2) = \lim_{k \rt \infty}d^{\tilde
C}_{\3}(z_1,z_2).$$ Hence it suffices to show the stability of the
Carath\'eodory distance.

 As before let $\Z \in C \subset \subset\Sigma_z$. for large
$k$, $\Z \in \3$. Let $\phi_k: \3 \rt \D$ be holomorphic maps such
that $\phi_k(z_1)=0$ and $$d^C_{\3}(\Z)=\rho(0, \phi_k(z_2)).$$

The family $\{ \phi_k \}$ is uniformly bounded above and since $
\3 \rt \Sigma_z$, all mappings $\phi_k$, $k \ge k_0$, are defined
on the compact set $C$. Thus there is a subsequence which will
still be denoted by $\phi_k$ so that  $\phi_k \rt \phi: \Sigma_z
\rt \overline \D$ and $ \phi (z_1)=0$. If $\vert \phi (z_0) \vert
=1$ for some $z_0 \in \Sigma_z$, then $\vert \phi (z) \vert
\equiv 1$ by the maximum principle. Thus $\phi: \Sigma_z \rt \D$
and in particular $\rho (0, \phi_k(z_2)) \rt \rho(0, \phi(z_2)).$
Therefore $d^C_{\3}(\Z) \le \rho (0, \phi(z_2)) + \ep \le
d^C_{\Sigma_z}(\Z)+ \ep$, which shows that
\begin{equation}\label{lll}
\limsup_{ \ k \rt \infty} d^C_{\3}(\Z) \le d^C_{\Sigma_z}(\Z)
\end{equation}
Conversely, working with the same subsequence that was extracted
above, we have:
$$ d^C_{\Sigma_z}(\Z) \le d^K_{\Sigma_z}(\Z) \le d^C_{\3}(\Z) +
\ep $$ for $k$ large. But
$$ d^K_{\3}(\Z)= d^K_{\1}( H_k^{-1}(z_1), H_k^{-1}
(z_2))$$ and since $H_k^{-1}(z_1), \ H_k^{-1} (z_2)$ are both
close to $p$, for $k$ large, it follows that
$$d^K_{\1}(H_k^{-1} (z_1), H_k^{-1}
(z_2)) \le d^K_{B_p(\eta_1) \cap \1}(H_k^{-1}  (z_1),
H_k^{-1}(z_2)) + \ep$$ where $B_p(\eta_1)$ is a small enough
neighborhood of $p$. Since $B_p(\eta_1) \cap \1$ is convex
Lempert's work shows that the Kobayashi and Carth\'eodory
distances coincide. Combining the aforementioned observation, we
get
\begin{equation}\label{ddd}
d^C_{\Sigma_z}( \Z) \le d^C_{B_p(\eta_1) \cap \1}(H_k^{-1}(z_1),
H_k^{-1} (z_2)) + 2 \ep
\end{equation}
To conclude, it is known (see ~\cite{pin2}) that the
Carath\'eodory distance can be localised near strongly
pseudoconvex points , exactly like the Kobayashi distance. hence
\begin{align}
d^C_{B_p(\eta_1) \cap \1}(H_k^{-1} (z_1), H_k^{-1}  (z_2))  & \le
(1+ \ep) \ d^C_{\1}(H_k^{-1} (z_1), H_k^{-1}(z_2)) \\ \notag & =
(1+ \ep) \ d^C_{\3} (\Z). \notag
\end{align}
 With this (\ref{ddd}) becomes
$$d^C_{\Sigma_z}( \Z) \le (1+ \ep) d^C_{\3}( \Z) +2 \ep.$$
Since $d^C_{\3}( \Z)$ are uniformly bounded by (\ref{lll}) it
follows that
\begin{equation}\label{kkk}
d^C_{\Sigma_z}( \Z) \le d^C_{\3}( \Z) +C \ep
\end{equation}
combining (\ref{lll}) and (\ref{kkk}), we see that
$$\lim_{k \rt \infty}d^C_{\3}( \Z) =d^C_{\Sigma_z}(\Z).$$
Hence the claim made in (\ref{ag1}) also holds for the
Carath\'eodory metric. The concluding arguments remain the same
in this case as well. This completes the proof of the proposition.
\end{proof}

Since the Kobayashi and Caratheodory distances coincide with a
constant multiple of the Bergman metric on $\Sigma_z$ and
$\Sigma_w$, it follows from \cite{gk1} that the limit map $\Phi :
\Sigma_z \rightarrow \Sigma_w$ is (anti)-biholomorphic.

\section{The boundary map is CR/anti-CR}
We prove Theorem \ref{mai} in this section. Throughout, $f: \1
\rt \2$ will denote a $C^1$-isometry of the Kobayashi or
Cartath\'eodory metrics which has a $C^1$-extension to $
\overline \Omega_1$.

Fix $p \in
\partial \1$. For the rest of this section we assume that
$p=f(p)=0$ and that the real normals to $\1$ and $\2$ at $p$ and
$f(p)$ are given by $\{ \tz= {\rm Im} \ z_n =0 \}$ and $\{ \tw =
{\rm Im} \ w_n =0 \}$. This can be achieved by composing $f$ with
transformations of the type in Lemma \ref{vit}.

Fix a sequence $\de_k \rt 0$ and define $x_k \in \1$ by $x_k =
(\tilde 0, - \de_k)$. Then $x_k \rt p$ as $k \rt \infty$. Because
of our choice of $x_k$, in the notation of Section 3, the map $h_k
=id$.

Recall that $f_k=g_k \circ f $ and $\Phi_k = T_k \circ f_k  \circ
T_k^{-1}$. More explicitly
$$\Phi_k(z)=\Biggl ( \frac {1}{\sqrt \dk}\tilde f_k(\sqrt \dk  \tz, \dk
\bz), \frac {1}{\dk}  {(f_k)_n}(\sqrt \dk \tz, \dk \bz) \Biggl)
$$

\begin{lem}\label{bor}
$\vert \tilde f_k(0) \vert =O(\de_k)$ {\rm as} $k \rt \infty$.
\end{lem}

\begin{proof}
Let $m$ be an upper bound for $\vert df \vert$ on $\overline
\Omega_1$. Now, noting that $\tilde f_k(x_k)=0$ and $\vert \tv
\vert \le \vert v \vert$ for any $v \in \Cn$, we have
\begin{align}\notag
\vert \tilde f_k( 0) \vert & = \ \vert \tilde f_k(x_k) - \tilde
f_k( 0) \vert \ \le \ \vert f_k(x_k) - f_k( 0) \vert \ \\ \notag
& \le C \ \vert f(x_k) - f( 0) \vert \le \ Cm \vert x_k  \vert \
\\ \notag & = \ C m \de_k  \notag
\end{align}
since $x_k=(0, -\de_k)$. In the second inequality we have used
the fact that $g_k \rt id$ in $C^ \infty$ on compact subsets of
$\Cn$.
\end{proof}

\vspace{2mm} \noindent The next lemma provides the crucial link
between the limit of the rescaled isometries and the derivative
of the boundary map. We clarify the notation used in the
statement and proof: First, even when we use complex notation,
all quantities will be regarded as entities on real Euclidean
space. In particular $\Cn$ is identified with ${\mathbb R}^{2n}$
by $z=(z_1,..,z_n)=(x_1+ \sqrt {-1}x_2,..,x_{2n-1}+ \sqrt
{-1}x_{2n}) \equiv (x_1,..x_{2n})$. Second, by the normalizations
made at the beginning of this section we note that, at $p$, the
decomposition $T_p \Omega_1= \Cn = H_p(\po_1) \oplus
H_p(\po_1)^{\perp}$ coincides with $\Cn= {\C}^{n-1} \oplus \C$.
Hence, by an abuse of notation, for $v \in T_p \Omega_1$, if
$v=v_H+v_N$ then $v=(v_H,v_N)$. Similar remarks hold for $f(p)$
and $\Omega_2$.

Next, by
 Proposition \ref{con}, a subsequence of $ \{ \Phi_k \}$ converges
to a (anti)-holomorphic automorphism $\Phi : \Sigma_z \rt
\Sigma_w$. For the statement of the lemma it helps to regard
$\Sigma_z$ and $\Sigma_w$ as subsets of $T_p \Omega_1$ and
$T_{f(p)} \Omega_2$ respectively.
\begin{lem}\label{bd}
With notation as above, for any $z=(\tz, \bz) \in \Sigma_z$, we
have $\widetilde {df_p}( \tilde z, 0) =d \tilde f_p( \tilde z,0)=
\tilde \Phi (z)$.
\end{lem}
\begin{proof}
The first equality is clear from the definitions. As for the
second, consider a map $ r: \Cn \rt \C ^{n-1} $ with $r=(r_1,
\ldots , r_{n-1})$. Given $ \de $ we write
$$\frac {1}{\sqrt \de} \ r(\sqrt
\de \tz, \de \bz) = \frac {r(\sqrt \de \tz, \de \bz) - r(\sqrt \de
\tz, 0)} {\sqrt \de}+ \frac {r(\sqrt \de \tz, 0) -r(0,0)} {\sqrt
\de}+ \frac {r(0,0)} {\sqrt \de}$$ By using the mean value
theorem for one-variable functions repeatedly, we can rewrite the
above equation as
$$ \frac {1}{\sqrt \de} \ r(\sqrt \de
\tz, \de \bz) = M \sqrt \de \bz  +N \tz + \frac {r(0,0)} {\sqrt
\de}.$$
 Here $M$ and $N$ are real matrices of sizes $(2n-2) \times 2$
and $(2n-2) \times (2n-2)$ respectively with entries $M_{lm}=
\frac {\partial r_l}{\partial x_m}(\eta_{lm}(\de))$ and $N_{ij}=
\frac {\partial r_i}{\partial x_j}(\xi_{ij}(\de))$. Also $\tz$ and
$z_n$ are regarded as column vectors of sizes $(2n-2) \times 1$
and $2 \times 1$ respectively.

The entries of $\eta_{ij}(\de) \in {\mathbb R}^{2n}$ lie between
the corresponding entries of $( \sqrt \de \tz , 0)$ and $(\sqrt
\de \tz, \de \bz)$. Similarly, the entries of $\xi_{ij}(\de)$ lie
between the entries of $(0,0)$ and $( \sqrt \de \tz , 0)$.

 Now apply this to $r= \tilde f_k$
and $\de=\de_k$ and let $k \rt \infty$. The first term goes to
zero since the entries of $M$ are bounded and the last term goes
to zero by Lemma \ref{bor}. Note that since $\{ g_k \}$ converges
to the identity map as $k \rt \infty$, we have $ \frac {\partial
(\tilde f_k)_i}{\partial x_j}(\xi^k_{ij}(\de)) \rt \frac
{\partial \tilde f_i}{\partial x_j}(0)$ by the continuity
assumption of $df$ on $\overline \Omega_1$. Hence the middle term
converges to $d \tilde f_p( \tilde z,0)$.

To complete the proof, we observe that since $\Phi_k \rt \Phi$,
$\frac {1}{\sqrt \dk}\tilde f_k(\sqrt \dk  \tz, \dk \bz) \rt
\tilde \Phi(z)$.

\end{proof}


We proceed with the proof of Theorem \ref{mai}. Recall that if
$\1$ and $\2$ are domains with smooth boundaries in $\Cn$, a
$C^1$ map $\phi: \po_1 \rt \po_2$ is said to be {\it CR} if, for
every $p \in \po_1$, the following two conditions are satisfied
$$d \phi_p(H_p(\po_1)) \subset H_{\phi (p)}(\po_2)$$
and
$$ d \phi_p \circ J_1 =J_2 \circ d \phi_p$$ where $J_1$ and $J_2$ are
the almost complex structures on $H_p(\po_1)$ and
$H_{\phi(p)}(\po_2)$. Similarly, an {\it anti-CR} map satisfies
$d \phi_p(H_p(\po_1)) \subset H_{\phi (p)}(\po_2)$ and $ d \phi_p
\circ J_1 = -J_2 \circ d \phi_p$ for every $p \in \po_1$.

In our case, $df$ satisfies the first condition by Lemma
\ref{hor}. We claim the second condition is satisfied due to
Lemma \ref{bd}. It follows from this lemma that the map $T:
{\mathbb B}^{n-1} \rt \C^{n-1}$ given by $T(\tz)= \tilde
\Phi(\tz,-1)$ is the restriction of the $\R$-linear map $
\widetilde {df_p}: \C^{n-1} \rt \C^{n-1}$.
On the other hand, $\tilde \Phi$ is holomorphic or
anti-holomorphic (since $\Phi$ is so). Combining these two
observations, it follows that $T$ is actually the restriction of a
$\C$-linear map. Hence,
\begin{equation}\label{lin}
\tilde \Phi (J_1(v),-1) = \pm  J_2 \tilde \Phi ( v,-1),
\end{equation}
for any $v \in {\mathbb B}^{n-1} $ and where $J_1, J_2$ denote
the almost-complex structures on $H_p(\po_1)$ and
$H_{\Phi(p)}(\po_2)$ respectively (note that we have used the
identification of the horizontal subspaces with $\C^{n-1}$).

This implies that $df_p$ is actually a $\C$-linear or conjugate
linear map on $\C^{n-1}=H_p(\po_1)$. More explicitly,
let $v \in H_p (\partial \1 )$. By scaling $v$ by some constant
$\alpha >0$, we can assume that $(\alpha v,-1) \in \Sigma_z$. By
Lemmas \ref{bd} and \ref{hor}, we have
$$ df_p(\alpha v,0)= \bigl(\tilde \Phi(\alpha v,-1 ),0 \bigr).$$
Using the $\C$-linearity/conjugate-linearity of $\Phi$ as in
(\ref{lin}), we have
$$df_p \bigl( J_1 (\alpha v),0 \bigr )= \bigl ( \tilde \Phi
 (J_1(\alpha v),-1), 0 \bigr )= \pm \bigl ( J_2 \tilde \Phi
( \alpha v,-1),0 \bigr)= \pm J_2 df_p( \alpha v,0).$$

Hence we conclude that the boundary map is CR/anti-CR.

Now we prove that $df_p:T_p \po_1 \rt T_{f(P)} \po_2$ is an
isomorphism. First, note that $d f_p \vert_{H_p(\po_1)}:
H_p(\po_1) \rt H_{f(p)}(\po_2)$ is invertible. To see this, let
${\mathbb H}_z := \Sigma_z \cap \{ (0,z_n): z_n \in \C \} $. Then
$\Phi(\Hy _z) \subset \Hy _w$ by Lemma \ref{bd}. On the other
hand, it can be checked that the induced Riemannian metrics on
$\Hy _z$ and $\Hy _w$ are just the hyperbolic metrics. From the
completeness of these metrics it follows that $\Phi(\Hy _z)= \Hy
_w$. But if $df_p$ is not injective on $H_p(\po_1)$, then there
would be a $(v,0) \in H_p(\po_1)$ such that $df_p(v,0)= 0$. By
scaling $v$ we can assume that $(v,-1) \in \Sigma_z$ and use
Lemma \ref{bd} to conclude that $\Phi(v,-1) \in \Hy _w$. This
contradicts $\Phi^{-1}(\Hy _w)=\Hy _z$.


Next, as in the proof of Lemma \ref{hor}, Equation $(\ref{ma})$
shows that $df_p(v)_N \neq 0$ for any $v \in H_p (\po_1) ^
\perp$. Hence $df_p : T_p \po_1 \rt T_{f(p)} \po_2$ is invertible
and $f$ is a CR/anti-CR diffeomorphism.

To conclude that $\1$ and $\2$ are biholomorphic we proceed as
follows: Note that the connectedness of $\po_1$ implies that $f:
\po_1 \rt \po_2$ is either CR or anti-CR everywhere. Let us
assume that $f$ is CR everywhere, the other case being exactly
similar. By ~\cite{pin1} it follows that there is a neighborhood
$U_1$ of $\po_1$ and a holomorphic mapping $F: U_1 \cap \1 \rt
\2$ such that $F$ is $C^1$-smooth upto $\po_1$ and $F=f$ on
$\po_1$. By Hartogs' theorem, $F$ extends to a holomorphic
mapping $F: \1 \rt \2$. Similarly $f^{-1}$ has a holomorphic
extension, say $G: \2 \rt \1$, which agrees with $f^{-1}$ on
$\po_2$. Since $f \circ f^{-1}=F \circ G =id$ on $ \po_2$, the
uniqueness theorem of ~\cite{pin1} forces $F \circ G =id$ on $\2$
and likewise $G \circ F =id$ on $\1$. Thus $\1$ and $\2$ are
biholomorphic. \hfill $\square$



\end{document}